\documentstyle[twoside]{article}

 \newcommand{\zr}[1]{\mbox{\hspace*{#1em}}}

 \newcommand{\CC}{\mbox{\zr{0.1}\rule{0.04em}{1.6ex}\zr{-0.30}{\sf C}}}
 \newcommand{\ZZ}{\mbox{\sf Z\zr{-0.45}Z}}
 \newcommand{\ZZs}{\mbox{\footnotesize \sf Z\zr{-0.45}Z}}
 
 \newcommand{\Id}{\mbox{{\bf 1}$\hspace{0.05em}\!\!$I}}
 
\input wigsym2.sty
\runninghead{
 M.~Baake and J.A.G.~Roberts  
}{
 Symmetries and reversing symmetries of trace maps
}
\titleline{
 $\mbox{SYMMETRIES AND REVERSING SYMMETRIES OF TRACE MAPS}$
}
\writers{ \sc Michael Baake }
\place{Institut f\"ur Theoretische Physik, Universit\"at T\"ubingen, \\
Auf der Morgenstelle 14, 72076 T\"ubingen, Germany}
\vglue 10pt
\centerline{\eightrm and} \vspace*{-0.1 truein}
\writers{\sc John A. $\!$G.~Roberts}
\place{Department of Mathematics, The University of Melbourne, \\
                 Parkville, Victoria 3052, Australia}
\vspace{0.35 truein}
\abstracts{\noindent
A (discrete) dynamical system may have various symmetries and 
reversing symmetries, which together
form its so-called reversing symmetry group.
We study the set of 3D trace maps 
(obtained from two-letter substitution rules)
which preserve the Fricke-$\!$Vogt invariant 
$\mbox{\small $I(x,y,z)$}$. 
This set of dynamical systems forms
a group $\mbox{\small $\cal G$}$ isomorphic with the 
projective linear (or modular) 
group $\mbox{\small $PGl(2,\ZZs)$}$. 
For such trace maps, we give a  complete characterization 
of the reversing symmetry group as a subgroup of 
the group $\mbox{\small $\cal A$}$
of all polynomial mappings that preserve $\mbox{\small $I(x,y,z)$}$.
}{}{}

\vspace{0.2 truein}


\section{Introduction}

Trace maps of two-letter substitution rules proved useful for the determination
of physical properties of various 1D systems$^{\ref{Kohmoto}}$ 
on non-periodic self-similar
structures such as the Fibonacci chain. 
These mappings found a systematic algebraic description
in terms of the free group with two generators and its automorphism 
(resp.\ homomorphism) group$^{\ref{Magnus},\ref{BGJ}}$; 
for a summary of theory and applications, see  Ref.~3.

Here, we are interested in trace maps as 3D dynamical systems, in particular in
those of the Nielsen class which preserve the Fricke-$\!$Vogt 
invariant$^{\ref{BGJ},\ref{Roberts}}$,
\vspace{-1mm} \begin{equation} \label{1.0} 
    I(x,y,z) \, = \,  x^2 + y^2 + z^2 - 2 x y z - 1.
\vspace{-1mm} \end{equation} 
The set ${\cal G}$ of {\em Nielsen trace maps} is a 
group$^{\ref{BGJ},\ref{Roberts}}$
isomorphic with the projective linear (or modular) group 
$PGl(2,\ZZ) \simeq Gl(2,\ZZ)/\{\pm \Id\}$.
In the standard representation$^{\ref{Magnus},\ref{BGJ}}$, 
every $F \in {\cal G}$ is an invertible polynomial
mapping of 3-space with integer coefficients. 
The group ${\cal G}$ can be characterized as
\vspace{-1mm}  \begin{equation} \label{2.0}
   {\cal G} \; = \;  \{ F \in \ZZ[x,y,z]^3 \,\,\, | \,\,\,
                 I(F(x,y,z)) = I(x,y,z) \,\, \mbox{and} 
                 \,\, F(1,1,1) = (1,1,1) \} \, ,
\vspace{-1mm} \end{equation}
with composition of mappings as group multiplication and 
${\cal G} \simeq PGl(2,\ZZ)$.

\section{Symmetries and reversing symmetries}

A trace map $F$ may possess certain {\em symmetries} (i.e., there may be
invertible transformations $S$ with $F\circ S = S \circ F$) which form a group,
${\cal S}(F)$. Typically, one specifies a subclass ${\cal B}$ of mappings
in which to look for possible symmetries, which we indicate here by
the notation ${\cal S}_{\cal B}(F)$. The choice of this
subclass is guided by key properties of the original
transformations, i.e., in our case, the trace maps. We shall give
a systematic description of ${\cal S}_{\cal A}(F)$, the group of
symmetries of Nielsen trace maps within
the set ${\cal A}$ of polynomial mappings that leave $I(x,y,z)$ invariant,
\vspace{-1mm} \begin{equation} \label{3.0}
   {\cal A} \; =  \;  \{ A \in \CC[x,y,z]^3 \,\,\, | \,\,\,
                 I(A(x,y,z)) = I(x,y,z) \}.
\vspace{-1mm} \end{equation}
It turns out that ${\cal A}$ is a {\em group} (hence, $A \in {\cal A}$ implies
$A$ invertible) which contains ${\cal G}$ as
a subgroup via the semidirect product
\vspace{-1mm} \begin{equation} \label{4.0}
   {\cal A} \; =  \; \Sigma \otimes_s {\cal G} \; \simeq \;
              \Sigma \otimes_s PGl(2,\ZZ) \, ,
\vspace{-1mm} \end{equation} 
where $\Sigma = \{Id,\sigma_1,\sigma_2,\sigma_3 \} \simeq C_2 \otimes C_2$ 
(Klein's 4-group)
is the normal subgroup$^{\ref{Magnus},\ref{Roberts}}$.
(Here and in what follows, we use $C_n$ for the cyclic group
of order $n$.) The involutions $\sigma_i$ are given by
$\sigma_1(x,y,z) = (x,-y,-z)$ etc.
The set ${\cal A}$ is quite remarkable: for any 
mapping $A \in {\cal A}$,
all polynomial coefficients are integers$^{\ref{Magnus},\ref{Roberts}}$,
wherefore we automatically have $A \in \ZZ[x,y,z]^3$.

{}Furthermore, $F \in {\cal G}$ can be {\em (weakly) reversible} (i.e., there is a 
reversing symmetry $G$ with $G\circ F\circ G^{-1}=F^{-1}$) or
{\em (strictly) reversible} (if $G^2=Id$, the identity mapping), see Ref.~4.
We will drop the historical
attributes ``weakly'' and ``strictly'', but distinguish trivial reversibility
(where $F^2 = Id$) from {\em true} reversibility (where $F$ is not an involution).
Reversibility is an important tool for the analysis of dynamical 
systems$^{\ref{RQ},\ref{Lamb}}$. 
The symmetries and reversing symmetries together form the so-called
reversing symmetry group$^6$, ${\cal R}(F)$. If $F$ with 
$F^2 \neq Id$ is reversible, ${\cal R}(F)$
is a {\em group extension} of the symmetry group 
${\cal S}(F)$ of index 2. Alternatively,
it can be seen as a group with a natural $\ZZ_2$-grading$^7$.
Here, we shall describe ${\cal R}_{\cal A}(F)$. 
This is not really a restriction --
see Ref.~4 for a first discussion of extensions, in particular to {\em 
arbitrary} polynomial mappings of finite order.

\section{Results}

Let us now formulate our results. We do this in two steps, giving the 
structure of the (reversing) symmetry group first within
${\cal G} \simeq PGl(2,\ZZ)$ and then within ${\cal A}$.
\vspace{1mm}
\begin{theorem} \label{t1}
 The symmetry group ${\cal S}_{\cal G}(F)$ for $F \in {\cal G}$ is: \\
  (1) ${\cal G}$ if and only if $F = Id$, \\
  (2) $C_2 \otimes C_2$ if and only if $F^2 = Id \neq F$, \\
  (3) $C_3$ if and only if $F^3 = Id \neq F$, \\
  (4) $C_{\infty}$ in all remaining cases, i.e., 
      whenever $F$ is not of finite order.
\end{theorem} 
\vspace{1mm}
The symbol $C_{\infty}$ is used for the infinite cyclic group generated from
one generator, i.e.~$C_{\infty} \simeq \ZZ$ as groups.
Unfortunately, we can only give a hint of the proof here, 
for details see Ref.~7. One has to
use ${\cal G} \simeq PGl(2,\ZZ)$ and to distinguish the elements of finite order
(possible is order 1, 2, or 3, where the proof is almost straightforward) 
from those of infinite order. In the latter case,
${\cal S}_{\cal G}(F)$ does certainly contain 
the subgroup $\{F^n \, | \, n \in \ZZ \}$, but often
more: $F$ may have roots within ${\cal G}$ which still give rise
to the group $C_{\infty}$. To show that
nothing else can happen, one uses the Dirichlet unit theorem of algebraic
number theory for quadratic fields. 
This gives the answer for all truly hyperbolic
or elliptic $2 \! \times \! 2$-matrices $M_F$, to be taken mod $\pm \Id$.
The remaining case, $\mbox{tr}(M_F) = 2$ with $\det(M_F)=1$, 
can be checked explicitly.
\vspace{1mm}
\begin{theorem} \label{t2}
 The reversing symmetry group ${\cal R}_{\cal G}(F)$ of $F \in {\cal G}$ is: \\
  (1) $C_{\infty}$ if and only if $F$ is not reversible, \\ 
  (2) ${\cal G}$ if and only if $F=Id$, \\
  (3) $D_2 \simeq C_2 \otimes C_2$ if and only if $F^2=Id \neq F$, \\
  (4) $D_3$ if and only if $F^3=Id \neq F$, or \\
  (5) $D_{\infty}$, if $F$ is truly reversible but not of finite order.
\end{theorem}  
\vspace{1mm}
Here, $D_n \simeq C_n \otimes_s C_2$ is the dihedral group, 
with $D_{\infty}$ being the extension to $n=\infty$. This Theorem
shows that all elements of finite order are reversible$^4$, which is trivial
only for order 1 and 2. If $F$ is an involution, we just have
${\cal S}_{\cal G}(F) = {\cal R}_{\cal G}(F)$.
The question whether a given element of infinite
order actually is reversible involves the representation of $\pm 1$ by integer 
quadratic forms and is, for any given $F$, decidable in finitely
many steps$^{\ref{BR}}$.

Now, the extension of
${\cal S}_{\cal G}(F)$ (${\cal R}_{\cal G}(F)$) to ${\cal S}_{\cal A}(F)$ 
(${\cal R}_{\cal A}(F)$), respectively, requires Eq.~(\ref{4.0}) and the
precise knowledge of how ${\cal G}$ acts on $\Sigma$.
With $\Sigma = \{ Id, \sigma_1, \sigma_2, \sigma_3 \}$,
$F \in {\cal G}$ acts -- via conjugation -- on $\Sigma$ as a
permutation $\pi_F^{} \in S_3$:
\vspace{-1.5mm} \begin{equation} \label{5.0}
   {}F \circ \sigma_i \circ F^{-1} \, = \, \sigma_{\pi_F^{}(i)}.
\vspace{-1.5mm} \end{equation} 
The permutation $\pi_F^{}$ is uniquely determined by $F$ and  easy
to calculate$^{\ref{Roberts},\ref{BR}}$. 
Such a permutation can only be of order
1, 2, or 3, so $(\pi_F^{})^6 = id$ (the identity in $S_3$)
for {\em all} $F \in {\cal G}$ and
\vspace{-1.5mm} \begin{equation}
   {\cal S}_{\cal A}(F) \, \subset \,  {\cal S}_{\cal A}(F^{\, 6}) \; = \;
    \Sigma \otimes_s {\cal S}_{\cal G}(F^{\, 6}) \; .
\vspace{-1.5mm} \end{equation} 
The exponent 6 can of course be replaced by the number 
$n=\mbox{ord}(\pi_F^{})$. Now, we have two possibilities: 
either ${\cal S}_{\cal G}(F^{\, 6}) = {\cal G}$ (if $F$ is of
finite order) or ${\cal S}_{\cal G}(F^{\, 6}) = {\cal S}_{\cal G}(F)
\simeq C_{\infty}$ (if $F$ is not of finite order). 
The latter statement is a consequence of the Cayley-Hamilton 
theorem applied to unimodular matrices. If $C_{\infty}$ is generated
by $F'$ (which could be a root of $F$ in ${\cal G}$) and if
$\pi^{}_{F'}=id$, the product in Eq.~(6) is direct.
 
Similarly, in case of {\em true} reversibility 
(i.e., $[{\cal R}_{\cal A}(F) : {\cal S}_{\cal A}(F)] = 2$), one obtains
\vspace{-1.5mm} \begin{equation} \label{6.0}
   {\cal R}_{\cal A}(F) \, \subset \, {\cal R}_{\cal A}(F^{\, 6}) \; = \;
    \Sigma \otimes_s {\cal R}_{\cal G}(F^{\,6}),
\vspace{-1.5mm} \end{equation}
where we can again simplify the right hand side
through the relation ${\cal R}_{\cal G}(F^{\, 6}) = 
{\cal R}_{\cal G}(F)$ if $F$ is not of finite order$^{4,7}$,
and through ${\cal R}_{\cal G}(F^{\, 6}) = {\cal G}$ if $F^{\, 6}=Id$.

However, one might not only be interested in the (reversing) symmetry
group of some power of $F$, cf.~Ref.~6, but in that of $F$ itself. 
To this end, we have to define
\vspace{-1.5mm} \begin{equation} \label{7.0}
   K^{}_{\Sigma}(F) \; := \; \{ g \in \Sigma \,\, | \,\, [g,F] = 0 \}
\vspace{-1.5mm} \end{equation}
which is a subgroup of $\Sigma$, hence one of the groups $\{Id\}$,
$\{Id,\sigma_1\}$, $\{Id,\sigma_2\}$, $\{Id,\sigma_3\}$, or $\Sigma$.
{}From the factorization property$^{\ref{Roberts}}$ 
of the semidirect product (\ref{4.0}), we then obtain 
\begin{theorem} \label{t3}
Let $F$ be a Nielsen trace map. The symmetry group within ${\cal A}$ is 

\vspace*{1mm}
\centerline{${\cal S}_{\cal A}(F) = K^{}_{\Sigma}(F) 
\otimes_s {\cal S}_{\cal G}(F)$}   

\vspace*{1mm}
\noindent while -- in case of true reversibility -- 
the reversing symmetry group is 

\vspace*{1mm}
\centerline{${\cal R}_{\cal A}(F) = K^{}_{\Sigma}(F) 
\otimes_s {\cal R}_{\cal G}(F)$.}
\end{theorem}
\vspace{1mm}
The products can be direct if the elements of
$K^{}_{\Sigma}(F)$ commute with all elements
of the symmetry group ${\cal S}_{\cal G}(F)$ resp.\
the reversing symmetry group ${\cal R}_{\cal G}(F)$
(which includes the case $K^{}_{\Sigma}(F) = \{Id\}$).
A comparison of Thm.~3 with Eqs.~(6) and (7) 
illustrates how the semidirect product structure
of ${\cal A}$ can lead to an element of ${\cal G}$ that has
fewer symmetries and reversing symmetries than certain of
its iterates. One example is provided by the Fibonacci
trace map$^4$ where $F^3$ commutes with every element
of $\Sigma$ but $F$ does not.

{}From the three Theorems, the structure of the (reversing) symmetry group
within ${\cal G}$ as well as within ${\cal A}$
is completely classified. Additionally, due to several examples,
we tend to believe that any $F$ which is not reversible within ${\cal A}$
is not reversible at all$^{\ref{Roberts},\ref{BR}}$.  

\section{Concluding remarks}

Let us briefly comment on general trace maps of two-letter substitution
rules. If they do not belong to the Nielsen class discussed above, they
cannot be of finite order and, even more, they are not globally invertible.
This also implies that non-Nielsen trace maps are never reversible in the
above sense. However, certain (generalized) symmetries and reversing
symmetries may show up in subspaces. This calls for a suitable generalization
of the (reversing) symmetry concept, in particular for the consideration
of non-invertible mappings, e.g., projectors in combination with 
an invertible transformation.
Though this gives a variety of interesting possibilities, we presently 
do not see anything close to a classification like that given above for the
Nielsen class.

\nonumsection{Acknowledgements}
It is a pleasure to thank J.S.W.~Lamb for interesting
discussions and communication of results prior to publication.
J.A.G.R.\ gratefully acknowledges the financial
support of the Australian Research Council through its
{}Fellowship Scheme.
M.B.\ would like to thank B.~Iochum and the CPT Luminy
for  hospitality and financial support during a stay in spring 1995
where this revised version of the manuscript has been completed.

\nonumsection{References}

{\small
\begin{enumerate}
\baselineskip=13pt
\parskip=1pt

\item \label{Kohmoto}
M.\ Kohmoto, L.P.\ Kadanoff and C.\ Tang, 
``Localization problem in one dimension: mapping and escape'',
Phys.\ Rev.\ Lett.\ {\bf 50} (1983) 1870--2; 
S.\ Ostlund, R.\ Pandit, D.~Rand, H.-J.\ Schellnhuber and E.D.\ Siggia, 
``One-dimensional Schr\"odinger equation with an almost 
periodic potential'', Phys.\ Rev.\ Lett.\ {\bf 50} (1983) 1873--6;
B.\ Sutherland, ``Simple system with quasiperiodic
dynamics: a spin in a magnetic field'',
Phys.\ Rev.\ Lett.\ {\bf 57} (1986) 770--3;
J.M.\ Luck, ``Frustration effects in quasicrystals:
an exactly soluble example in one dimension'',
J.\ Phys.\ {\bf A20} (1987) 1259--68;
J.M.\ Luck, H.\ Orland and U.\ Smilansky, 
``On the response of a two-level quantum system to a 
class of time-dependent quasiperiodic perturbations'',
J.\ Stat.\ Phys.\ {\bf 53} (1988) 551--64.  

\item \label{Magnus}
R.D.\ Horowitz, ``Induced automorphisms on Fricke
characters of free groups'',
Trans.\ Am.\ Math.\ Soc.\ {\bf 208} (1975) 41--50; 
W.\ Magnus, ``Rings of Fricke characters and automorphism
groups of free groups'', Math.\ Z.\ {\bf 170} (1980) 91--103; 
J.-P.\ Allouche and J.\ Peyri\`ere, 
``Sur une formule de r\'ecurrence sur les traces de produit
de matrices associ\'es \`a certaines substitutions'',
C.\ R.\ Acad.\ Sci.\ Paris {\bf 302}(II) (1986) 1135--6; 
J.\ Peyri\`ere, ``On the trace map for products of matrices
associated with substitutive sequences'',
J.\ Stat.\ Phys.\ {\bf 62} (1991) 411--4; 
P.~Kramer, ``Algebraic structures for one-dimensional
quasiperiodic systems'', J.\ Phys.\ {\bf A26} (1993) 213--28;
P.~Kramer, ``Fricke-Klein geometry for the group $Sl(2,\CC)$'',
J.\ Phys.\ {\bf A26} (1993) L245--50;
J.\ Peyri\`ere, Wen Zhi-Xiong and Wen Zhi-Ying, 
``Algebraic properties of trace mappings associated with
substitutive sequences'', Modern Mathem.\ (China) (1993),
in press; P.\ Kramer and J.\ Garcia-Escudero,
``Automorphisms of free groups and quasicrystals'',
this volume.

\item \label{BGJ}
M.\ Baake, U.\ Grimm and D.\ Joseph, ``Trace maps,
invariants, and some of their applications'', Int.\ J.\ Mod.\
Phys.\ {\bf B7} (1993) 1527--50.

\item \label{Roberts}
J.A.G.\ Roberts and M.\ Baake, ``Trace maps as 3D reversible
dynamical systems with an invariant'', 
J.\ Stat.\ Phys.\ {\bf 74} (1994) 829--88;
J.A.G.\ Roberts and M.\ Baake, ``The dynamics of trace maps'', 
in: {\em Hamiltonian Mechanics: Integrability and Chaotic Behaviour}, 
ed.\ J.\ Seimenis, NATO ASI
Series B: Physics, (Plenum Press, New York, 1994) 275--85.

\item \label{RQ}
J.A.G.\ Roberts and G.R.W.\ Quispel, ``Chaos and
time-reversal symmetry. Order and chaos in reversible
dynamical systems'', Phys.\ Rep.\ {\bf 216} (1992) 63--177.

\item \label{Lamb}
J.S.W.\ Lamb, ``Reversing symmetries in dynamical systems'',
J.\ Phys.\ {\bf A25} (1992) 925--37;
J.S.W.\ Lamb and G.R.W.\ Quispel, ``Reversing k-symmetries in
dynamical systems'', Physica D {\bf 73} (1994) 277--304.

\item \label{BR}
M.\ Baake and J.A.G.\ Roberts, ``Reversing symmetry group
of $Gl(2,\ZZs)$ and $PGl(2,\ZZs)$ matrices with connections 
to cat maps and trace maps'', J.\ Phys.\ {\bf A30} (1997) 1549--73.

\end{enumerate} }\hspace*{3cm}
\end{document}